\documentclass{amsart}

\usepackage{amsmath}%
\usepackage{amsfonts}%
\usepackage{amssymb}%
\usepackage{amsthm}
\usepackage{graphicx}
\usepackage{tikz-cd}
\usepackage{hyperref}
\usepackage[T1]{fontenc}
\usepackage{xcolor}[table]

\usepackage{geometry}
\DeclareUnicodeCharacter{3000}{ }

\theoremstyle{plain}

\newtheorem*{acknowledgement*}{Acknowledgement}

\newtheorem{definition}{Definition}

\numberwithin{equation}{section}
\newtheorem{theorem}{Theorem}[section]
\newtheorem{lemma}[theorem]{Lemma}
\newtheorem{proposition}[theorem]{Proposition}
\newtheorem{corollary}[theorem]{Corollary}

\theoremstyle{remark}
\newtheorem{remark}[theorem]{Remark}
\newtheorem{example}[theorem]{Example}
\numberwithin{equation}{section}

\DeclareMathOperator{\tr}{Tr}
\DeclareMathOperator{\rtr}{\mathcal{R}eTr}
\DeclareMathOperator{\pr}{Pr}
\DeclareMathOperator{\prv}{Pr_{v}}
\DeclareMathOperator{\prh}{Pr_{h}}
\DeclareMathOperator{\diag}{diag}

\DeclareMathOperator{\grad}{grad}

\DeclareMathOperator{\ver}{Ver}
\DeclareMathOperator{\hor}{Hor}

\newcommand{\Mn}{\mathcal{M}(n,\mathbb{R})}
\newcommand{\Mnm}{\mathcal{M}(n, m, \mathbb{F})}

\newcommand{\Mnmo}{\overset{\circ}{\mathcal{M}}(n, m, \mathbb{F})}

\DeclareMathOperator{\strat}{\textit{\dj}}

\newcommand{\herm}{\mathcal{H}(n)}
\newcommand{\hermo}{\overset{\circ}{\mathcal{H}}}
\newcommand{\algsu}{\mathfrak{su}(n)}

\newcommand{\ib}{\mathbf{i}}
\newcommand{\jb}{\mathbf{j}}
\newcommand{\kb}{\mathbf{k}}

\title[Eigenvalue as submersion]{Eigenvalue processes in light of Riemannian submersion and gradient flow of isospectral orbits}
\author{Huang Ching-Peng}
\email{cphuang@brown.edu}

\begin{document}
	\maketitle
	
	\begin{abstract}
		We prove eigenvalue processes from dynamical random matrix theory including Dyson Brownian motion, Wishart process, and Dynkin's Brownian motion of ellipsoids are results of projecting Brownian motion through Riemannian submersions induced by isometric action of compact Lie groups, whose orbits have nonzero mean curvature, which contributes to drift terms and is the log gradient of orbit volume function, showing in another way that eigenvalues collide whenever the fibre is degenerate. We thus provide a unified treatment and better connection between eigenvalue processes in different settings with the language of Riemannian geometry.
		
		Under such interpretation, we see how we can naturally recover eigenvector processes and derive $\beta$ process such as $\beta$-Dyson Brownian motion for general $\beta>0$.

		\textbf{KEYWORDS}: eigenvalue process, symmetric space, mean curvature flow, gradient flow, isospectral manifold, random matrix ensemble
		
	\end{abstract}

	\tableofcontents

	\section{Introduction}
	
	\subsubsection{Motivation}
	
	A major part of random matrix theory (RMT) is dedicated to behaviour of eigenvalues. Taking the dynamical point of view, Dyson \cite{dyson1962brownian} wrote down his famous stochastic process describing evolutions of eigenvalues of a ``hermitian Brownian motion'', $H_t$, whose entries are independent Wiener processes. The eigenvalue process is later known as the Dyson Brownian motion: 
	\begin{equation}\label{DBM}
		d\lambda_a = d\omega_a + \sum_{b\neq a} \frac{\beta dt}{\lambda_a - \lambda_b}.
	\end{equation}
	
	Similarly, the Wishart process, which comes from singular values of a rectangular matrix $M_t$ with IID Brownian motion entries was derived by Bru \cite{bru1989diffusions}, 
	
	\[ d\lambda_i = 2\sqrt{\lambda_i}d\omega_i + (\sum_{j\neq i}2\beta\frac{\lambda_i + \lambda_j}{\lambda_i - \lambda_j} + n)dt,\] where $\lambda$'s are eigenvalues of $M_tM_t^*$
	
	Dynkin's Brownian motion of ellipsoids, obtained as singular values of stochastic process compatible with the Lie structure of $GL(n)$, i.e.
	\[dG_t = G_t\strat W_t,\] where $W_t$ is a standard $n\times n$ Wiener process, and the eigenvalues of $G_tG_t^*$ was shown in \cite{article} by Norris and Rogers. The derivations of all three processes involve no more than elementary matrix operations and stochastic calculus. However, unified theory was not patently available despite the similarity among the three, such as the ``repulsive'' force in the drift term, e.g. $\frac{\beta dt}{\lambda_a - \lambda_b}$ in Dyson BM \ref{DBM}. A main purpose of this paper is to show how differential geometry is a tool to systematically describe eigenvalue processes.

	\subsubsection{Underlying geometry}
	
	Fortunately, suitable tools are available. The main idea is that:
	
	\begin{quotation}
		 Brownian motion gains a mean curvature flow term through Riemannian submersion. C.f. Theorem \ref{BM projection formula}
	\end{quotation}

	This is applicable once we realise projecting to eigenvalues or singular values, under proper setups, is a quotient by unitary groups, and in general quotient by a Lie group $K$ acting as isometry induces a Riemannian submersion. C.f. Theorem \ref{quotient submersion} 
	
	Our key result is summarised as the following:
	
	\begin{quotation}
		The eigenvalue processes including Dyson Brownian motion, Wishart process, and Dynkin's Brownian motion are all consequences of Brownian motion mapped through Riemannian submersion.
	\end{quotation} 
	
	In fact, the mean curvature is the log gradient of the ($K$-)orbit volume function. (C.f. \cite{pacini2003mean}.) (E.g. For Dyson BM, the volume function is proportional to the Vandermonde $\prod_{i \neq j}(\lambda_i-\lambda_j)$.) Not only does it simplify the calculation but interprets the repulsion as a force to increase the orbit volume, and we observe
	\begin{quotation}
		eigenvalues collide when the orbit of $K$ is degenerate, i.e. of lower dimension and hence has zero volume. 
	\end{quotation} 
	
	In Section \ref{EV section}, to see the eigenvector evolution, we study the local model of Riemannian submersion in Section \ref{semidirect BM} for a complete geometric treatment. 

	 Dyson examined the process over the field $\mathbb{R}, \mathbb{C}$, and the quaternion $\mathbb{H}$ respectively, resulting in what he referred to as ``three-fold'' for $\beta = 1, 2, 4$ in \ref{DBM}. Later, matrix models whose eigenvalues coincides with Dyson BM with general $\beta>0$ were constructed. (See e.g. \cite{unterberger2018global}.) Continuing what we have in \cite{huang2022dyson}, we see how the same construction works for Wishart and Dykin in a more general setting in Section \ref{general beta section}.
	
	\subsubsection{Structure of the paper}
	
	To begin, we introduce main tools that we use, followed by a list of notations and conventions. Before the main content, we provide by a warm-up example to familiarise the reader with some key elements.
	
	Next, we derive formulae for the eigenvalue processes, eigenvector processes, and construction for general $\beta$-processes completely in the Riemannian language, as outlined in the previous section. 
	
	In the appendix, we provide calculation for the mean curvature terms by definition.  
	
	We would like to note that, RMT is traditionally a subject heavily in analysis, and we focus on the interpretation of it in geometry. Careful analysis is abundant in the literature. See \cite{huang2022dyson} for an introduction to the nicety of the subject and reference.
	
	\subsubsection{Background knowledge}
	
	We assume basic knowledge of stochastic calculus and Riemannian geometry including Lie groups and Lie algebras. Some advanced terms might be used but should not obstruct the understanding if the reader has not encountered them. 
	
	For constructions of stochastic process on manifolds, we only mention in a conceptual manner. Readers can refer to standard text for the subject, e.g. \cite{hsu2002stochastic} and \cite{ikeda2014stochastic}.

	\subsection{Tools}
	
	The main tools we use are the following theorems.
	
	\begin{theorem}\label{BM projection formula}(\cite{je1990riemannian}\cite{watson1973manifold})
		Let $\phi: M \to N$ be a Riemmanian submersion with relative dimension, $k$. Let $\vec{H}$ denote the vector field on $M$ defined as the mean curvature of the fibres $\phi^{-1}(y)$ for $y \in N$. Suppose $\phi_*(H(x))$ is constant along each fibre, i.e. $\phi_*(H(x_1)) = \phi_*(H(x_2))$ if $x_1, x_2 \in \phi^{-1}(y)$; then if $W_M$ is a Brownian motion on $M$, \[d\phi(W_M) = dW_N - \frac{k}{2}\phi_*(H(W_M))dt,\] where $W_N$ is a Brownian motion on $N$. 
	\end{theorem} 
	
	In the case described in the theorem above, $H$ is a vector field on $M$, and we say it \textit{descends} to a vector field on $N$. 

	The submersion of interests in this paper are all quotients of compact Lie groups. Naturally, the it offers nice symmetries.
	
	\begin{proposition} \label{quotient submersion}(C.f.\cite{pacini2003mean})
		Let $M$ be a Riemannian manifold such that a compact Lie group $K$ acts as a subgroup of isometries of $M$. Then $M/K$ is equipped with a natural quotient metric such that $\phi: M \to M/K$ is a Riemannian submersion, and the orbit mean curvature field descends to a vector field on the quotient $M/K$.
	\end{proposition}
	
	The history of mean curvature flows dates back at least to the famous Plateau problem and variation of surface area with fixed boundary. In the group quotient scenario, we have a characterisation in this vein as follows.
	
	\begin{proposition} \label{log vol grad}(C.f. \cite{pacini2003mean} Proposition 1)
		Let $M$ be a Riemannian manifold such that a compact Lie group $K$ acts as a subgroup of isometries of $M$. The volume of orbit $vol(K\cdot m)$ is a smooth function on $M$, and the vector field of orbit mean curvatures is \begin{equation}
			\vec{H}(m) = -\grad\log vol(K\cdot m)
		\end{equation} for $m\in M$.
	\end{proposition}

	For group orbits, the volume function can be calculated as follows.
	
	\begin{lemma}\label{volume of orbit}(C.f. \cite{pacini2003mean} Proof of Proposition 1)
		Let $(M,g)$ be a Riemannian manifold such that a compact Lie group $K$ acts as a subgroup of isometries of $M$, and let $\{Z_\alpha\}_\alpha$ be a basis of $\mathfrak{k} = Lie(K)$, the volume of the orbit of a point $m\in M$ is \begin{equation}
			vol(K\cdot m) = \sqrt{\det i^*g}\int_K \wedge_\alpha Z_\alpha^*.
		\end{equation} Here $i^*g$ is the metric, in terms of the basis $\{Z_\alpha\}_\alpha$, pulled back from the inclusion $i: K\cdot m \hookrightarrow M$, and $Z_\alpha^*$ is the dual of $Z_\alpha$. 
	\end{lemma}

	The orbit volume function above is usually in terms of the volume of the homogeneous space $K/Stabilizer(m)$, which is constant for generic orbits. Taking the log gradient eliminates this constant, and we do not need the knowledge of this specific quantity to calculate the mean curvature. 

	\subsection{Notations and conventions} \label{notations}
\begin{itemize}
	\item $E_{ij}$ is the matrix with $1$ at the $(i,j)$-th position. $A_{ij} = E_{ij} - E_{ji}$, $S_{ij} = E_{ij} + E_{ji}$, and $D_i = E_ii$, sometimes up to a constant for convenience. We overload the operator $\diag$ so it takes a matrix to vector and vector to matrix, as commonly used in programming languages.
	
	\item $\mathcal{M}(n, \mathbb{F})$ or $\mathcal{M}_\mathbb{F}(n)$ is the set of $n\times n$ matrices over $\mathbb{F}$, or simply $\mathcal{M}(n)$ when the context is clear. Similarly $\mathcal{M}(n,m, \mathbb{F})$ is the set of $n\times m$ matrices. $\overset{\circ}{\mathcal{M}}$ means the subset consisting of element having distinct singular values.  
	
	\item $\mathcal{H}(n)$ denotes the set of $n\times n$ (complex) hermitian matrices and $\hermo$ the subset having distinct eigenvalues. We only used the notation for complex numbers although one could naturally write $\mathcal{H}(n, \mathbb{F})$ to include real numbers and the quaternion.
	
	\item For transpose of a real matrix or conjugate of a quaternion matrix $A$, we use the same notation for complex conjugate $A^*$ since the generalisation is natural in our narration.  
	
	\item Manifolds are assumed to be real, smooth, finite dimensional, possibly open or with boundary if not specified. 
	
	\item We adopt Einstein's summing convention.
		
	\item The unitary group $U(n)$ could be over field $\mathbb{F} = \mathbb{R}, \mathbb{C}$, or the quaternion $\mathbb{H}$, where $\ib, \jb, \kb$ are the imaginaries. The orthogonal group $O(n) = U(n, \mathbb{R})$ is always over real numbers.
	
	\item The differential of a smooth map, $\phi$, between smooth manifolds is written as $\phi_*$ instead of $d\phi$ to avoid confusion with stochastic differentials. 
	
	\item For stochastic processes, we often suppress the time variable $t$ for simplicity. For example, we write $X$ instead of $X_t$ when it is clear in the context. 
	
	\item For It\^o differentials we use $d$ while for Stratonovich differential we use $\strat$. The more common notation for the later is $\circ d$ and sometimes $\partial$ or $\delta$. We avoid $\circ d$ to prevent confusion with composition of maps.
	
	\item $\vec{H}(x)$ (sometimes without the vector arrow when not emphasizing it is a vector) denotes the mean curvature vector at a point $x$. We adopt the convention to define the mean curvature to be the trace of the second fundamental form divided by the dimension of the submanifold. Note that in most cases through out the paper, $\vec{H}$ is a vector field that descends to a vector field $\vec{J}(x)$ on the quotient manifold.  
	
	\item \textit{MCF} is short for \textit{mean curvature flow}, often scaled and reversed when it should cause confusions.
	
	\item\textit{BM} is short for \textit{Brownian motion}.
	
	\item O.n.b. is short for \textit{orthonormal basis}. 
	
	\item The spectrum of a matrix can refer to \textit{eigenvalues} or  \textit{singular values} when it is clear in the context. An \textit{isospectral set/orbit/manifold} is the subset of all matrices sharing the same spectrum.

	\item Given a Riemannian submersion $\phi: M\to N$, $\ver_xM$ and $\ver_xM$ denote the vetical and the horizontal spaces w.r.t. $\phi$ at $x\in M$, i.e. the tangent and the normal spaces to the fibre of $x$, respectively. 
	
	\item The \textit{Frobenius metric} of matrices over real, complex, or quaternion numbers is the standard Euclidean norm, conveniently defined as the real trace $\langle A,B\rangle_{Fr} = \rtr AB^*$.

\end{itemize}

\subsection{Warm-up example: normal metric on isospectral orbits}

	Before we start with the main content, we see a simpler example for Brownian motion through submersion that also introduces properties of isospectral orbits. (See also \cite{bloch1997double}\cite{https://doi.org/10.48550/arxiv.2205.06737}\cite{edelman1998geometry})

	\begin{example}(BM on flag manifold represented as isospectral orbit with the normal metric)
		
		Let \[\Lambda = \begin{bmatrix}
			\lambda_1 I_{n_1} &&\\
			& \ddots & \\
			&& \lambda_r I_{n_r}
		\end{bmatrix}, \] where $\lambda_k$'s are distinct real numbers and $\sum n_k = n$. The orbit of $\Lambda$ under the adjoint action of $O(n)$ is the manifold consisting of real symmetric matrices which diagonalise to $\Lambda$, \[\mathcal{O} = \{Q\Lambda Q^*| Q \in O(n)\}.\] 
	
		The canonical metric of $O(n)$ is the one obtained from the Killing form, which is equivalent to the submanifold metric from the usual embedding $O(n) \hookrightarrow \mathbb{R}^{n^2}$. 
		
		$\mathcal{O}$ is diffeomorphic to the flag manifold, $Flag(n_1, \cdots, n_r) = O(n)/O(n_1)\times\cdots\times O(n_r)$. The \textit{normal metric} on $\mathcal{O}$ is induced by the submersion $O(n) \to \mathcal{O}$. Since all fibres are of the same volume, they are totally geodesic by Theorem \ref*{log vol grad}, and consequently the Brownian motion on $O(n)$ projects to that on $\mathcal{O}$.
		
		The canonical BM on $O(n)$ can be written as the right-invariant process $Q$ such that
		\begin{equation}
			dQ = \strat A Q,
		\end{equation} where $A$ is skew-symmetric and has IID Wiener processes as entries, i.e. BM on the Euclidean space $\mathfrak{so}(n)$. Its projection to $\mathcal{O}$ is $H = Q\Lambda Q^*$, which satisfies
		\begin{align}
			dH & = \strat AQ\Lambda Q^* - Q\Lambda Q^*\strat A = [\strat A, H] \\
				& = [dA,H] + \frac{1}{2}[dA,dH].
		\end{align} For simplicity, let us only consider the case where $n_i = 1$ for all $i$, and $\mathcal{O}$ is then diffeomorphic to $O(n)$ quotient by its maximal torus $T$. 
	
		 From the above equations, 
		\begin{equation}
			dAdH = -(n-1)Hdt - dAHdA,
		\end{equation} and the $(i,l)$-th entry of the latter is
		\begin{align}
			da_{ij}h_jkda_{kl} & = (1-\delta_{il})da_{il}h_{li}da_{il} + \delta_{li}da_{ij}h_{jj}da_{ji} \\
				& = h_{li}dt - \delta_{li}\sum_{j\neq i}h_{jj}dt.
		\end{align}
		So \begin{align}
			dAHdA	& = (H - \diag(H))dt - ((\tr H)I -\diag(H))dt \\
			     	& = (H - (\tr H)I) dt.
		\end{align} Hence 
		\begin{equation}
			dAdH = -(n-1)Hdt - (H - (\tr H)I) dt = (-nH + (\tr H)I)dt,
		\end{equation} which is symmetric. Therefore 
		\begin{equation}
			[dA, dH] = dAdH - dHdA = dAdH + (dAdH)^* = 2dAdH.
		\end{equation} In conclusion,
		\begin{equation}
			dH = [dA,H] + dAdH = [dA,H] - (nH - (\tr H)I)dt.
		\end{equation} We have an autonomous equation for BM on the isospectral manifold. 
		
	\end{example}

	\section{Three eigenvalue processes}
	
	The derivation all three cases follow the same road map.
	\begin{enumerate}
		\item Given a matrix process $M_t$, we want to find the SDE for its eigenvalues.
		
		\item Identify the metric under which $M_t$ is a Brownian motion.
		
		\item Identify the group $K$ acting on $M_t$ such that each orbit is exactly an isospectral manifold.
		
		\item Verify the group $K$ acts isometrically and conclude the quotient induced by the action is a Riemannian submersion.
		
		\item Choose a basis of $K$ and pull back the metric from each orbit. 
		
		\item Compute the orbit volume function with the pull-back metric tensor. 
		
		\item  Compute the gradient of the log orbit volume to obtain the drift for the eigenvalues.
		
		\item Identify the quotient metric structure on the manifold of spectra and write down the SDE for BM on under it.
		
		\item The eigenvalue process is BM on the quotient space plus the drift term according to Theorem \ref{BM projection formula}.
	\end{enumerate}

	\subsection{Dyson Brownian motion}
	
	The space of $n\times n$ hermitian matrices, $\mathcal{H}(n, \mathbb{C})$ is a linear submanifold of $\mathcal{M}_\mathbb{C}(n) \cong \mathbb{C}^{n^2}$, equipped with the Frobenius (Euclidean) norm. We can take the field $\mathbb{F}$ to be $\mathbb{R}$, $\mathbb{C}$, or the quaternion $\mathbb{H}$, and by ``hermitian'', we include real symmetric, complex hermitian, and quaternion hermitian cases, denoted $\mathcal{H}(n,\mathbb{F})$, or simply $\mathcal{H}(n)$ when not specifying the field. We denote $(\hermo, Fr)$ as the open dense submanifold consisting of hermitian matrices with simple spectra. The tangent spaces may be identified with $\herm$.
	
 	Let
	\begin{equation}
		A_{ab} = \frac{E_{ab} - E_{ba}}{\sqrt{2}},
	\end{equation} 
	\begin{equation}
		S_{ab} = \frac{E_{ab} + E_{ba}}{\sqrt{2}},
	\end{equation} for $a\neq b$, and 
	\begin{equation}
		D_a = E_{aa}.
	\end{equation} BM on the Euclidean space $\herm$ is hence 
	\[W^{a}D_a + W^{ab}S_{ab}\] for $\mathbb{R}$,
	\[W^{a}D_a + W^{ab}S_{ab} + \ib\hat{W}^{cd}A_{cd},\] for $\mathbb{C}$, and 
	\[W^{a}D_a + W^{ab}S_{ab} + \ib\hat{W}^{cd}A_{cd} + \jb\hat{W}^{ef}A_{ef} + \kb\hat{W}^{gh}A_{gh}\] for $\mathbb{H}$, where the $W$'s are independent standard BM. 
	
	Consider the group action of $U(n)$ on $\herm$, $\rho: (U, S) \mapsto USU^*$ for $U\in U(n)$ and $S \in \hermo$. Since \begin{equation}
		\rtr AB^* = \rtr UAU^*UB^*U^*, 
	\end{equation} for $A$, $B \in T_S \hermo$ this is an isometry. Each orbit is exactly an isospectral orbit, those in $\herm$ having identical eigenvalues. We write the mapping \begin{equation}
	\rho^S: U(n) \to U(n)\cdot S.
	\end{equation} for the projection to an orbit.
	
	Next we consider the quotient map \begin{equation}
		\phi: \hermo \to \hermo/U(n)
	\end{equation} under such action and calculate the vertical and horizontal spaces w.r.t. the quotient. Since $U(n)$ acts by isometry, we may pick a diagonal element $\Lambda = \diag(\lambda_1, \cdots, \lambda_n)$ in the orbit to simplify the calculation.

	The Lie algebra $\mathfrak{su}(n)$ consists of traceless anti-hermitian matrices. Let $E_{ab}$ be the matrix with $1$ as the $(a,b)$-th entry and zero elsewhere. Then $\mathfrak{su}(n)$ is spanned over $\mathbb{R}$ with basis \[\mathcal{B} = \{\ib(D_a-D_b)\}\cup\{A_{ab}\},\] or
	\[\mathcal{B} =\{\ib(D_a-D_b)\}\cup\{A_{ab}\}\cup\{\ib S_{ab}\},\] or
	\[\mathcal{B} =\{\ib(D_a-D_b)\}\cup\{A_{ab}\}\cup\{\ib S_{ab}\}\cup\{\jb S_{ab}\}\cup\{\kb S_{ab}\}\] with $a<b$, for each field. 
	\vfil
	\noindent\textbf{Claim}: $\ver_\Lambda \hermo = \{[A,\Lambda]|A\in \mathfrak{su}(n) \text{ with } \diag(A) = 0\}$
	\begin{proof}
		For $A\in \algsu$, let $U_t = \exp(tA)$. Then 
		\[\rho^L_*(A) = \frac{d}{dt}|_{t=0}U_t\Lambda U_t^* = A\Lambda + \Lambda A^* = A\Lambda - \Lambda A.\]
		
		In particular, if $A$ is diagonal, $[A, \Lambda] = 0$.
	\end{proof}

	It is now obvious (real) diagonal matrices are orthogonal to the vertical space and have complementary dimensions. Hence we conclude \begin{equation}
		\hor_\Lambda \hermo = \{\text{real diagonal matrices}\}.
	\end{equation} The quotient manifold $\hermo/U(n)$ is parametrised by a single chart $\mathcal{W}= \{ (\lambda_1, \cdots, \lambda_n)|\lambda_1 > \cdots > \lambda_n\} \subset \mathbb{R}^n$ (i.e. $A_n$-Weyl chamber) under the mapping $\Lambda = \diag(\lambda_1, \cdots, \lambda_n) \mapsto U(n)\cdot \Lambda$; hence we easily check that the metric is Euclidean.

	We have now:
	\begin{proposition}
		The mapping $\phi: \hermo \to \hermo/U(n)$ is a Riemannian submersion where the quotient manifold is isometric to the $A_n$-Weyl chamber in $\mathbb{R}^n$. 
	\end{proposition}
	
	With the basis derived above, we are ready for the computation of the eigenvalue process, which is the image of Brownian motion through the quotient here.
	 
	\begin{proposition}\label{Dyson orbit volume}
		The orbit volume function of the action $U(n) \curvearrowright\hermo$ is \[c\cdot\prod_{a\neq b} (\lambda_a-\lambda_b)^\beta\] with $\beta = 1, 2, 4$ for real, complex, and quaternion respectively, where $c$ is a constant only dependent on $n$.
	\end{proposition}
	
	\begin{proof}
		$A_{ab}$ and $\{S_{ab}\}$ are a $U(n)$-invariant fields for the flag manifold $U(n)/T$, where $T$ is the maximal torus. Under $\phi$, \begin{equation}
			A_{ab} \mapsto (\lambda_b-\lambda_a)A_{ab}, \text{ and } \ib S_{ab} \mapsto \ib(\lambda_b-\lambda_a),
		\end{equation} and similar for $\jb$ and $\kb$. The pull-back metric only has diagonal terms $(\lambda_a-\lambda_b)^\beta$ with $\beta = 1, 2, 4$ for real, complex, and quaternion respectively. Hence the volume function is \begin{equation}
			vol(U(n)\cdot \Lambda) \propto \prod_{a\neq b}(\lambda_a - \lambda_b)^\beta
		\end{equation}
	\end{proof}
	
	Combining all, we recover the Dyson BM.
	
	\begin{corollary}\label{Dyson BM formula}
		The eigenvalue processes of the hermitian BM, i.e. Dyson BM, satisfy \[d\lambda_a = d\omega_a + \sum_{b\neq a} \frac{\beta dt}{\lambda_a - \lambda_b}, \] where $\omega_a$'s are independent standard BM.
	\end{corollary}
	
	\begin{proof}
		The quotient space has Euclidean metric and hence the form of BM. 
		
		The mean curvature is obtained by computing 
		\begin{equation}
			\frac{\partial}{\partial \lambda_a}\log \prod_{a\neq b}(\lambda_a - \lambda_b)^\beta = \sum_{b\neq a} \frac{\beta}{\lambda_a - \lambda_b}.
		\end{equation} 
	\end{proof}

\subsection{Wishart process}

The Wishart process starts with the $n\times m$-dimensional standard Brownian motion $W \in \mathbb{F}^{nm}$ for $n>m$, in other words an $n\times m$ matrix with independent entries of standard $\mathbb{F}$-BM. The eigenvalue process of interest is those of $WW^*$.  

We thus consider singular values of $W$. Let $K = U(n)\times U(m)$ act on $(\Mnm,Fr) \cong (\mathbb{F}^{nm}, Euc)$ by $(U, V): M \mapsto UMV^*$, and let $\rho^M: K \to K\cdot M \subset \mathbb{F}^{nm}$. Apparently, $K$ is compact and acts by isometry:
\[\rtr (UAV^*)(UBV^*)^* = \rtr AB^*\]
  for any $A$ and $B \in T_M\Mnm$. We consider the mapping $\phi: \Mnmo \to \Mnmo/K$, restricting our focus to the open locus where singular values are distinct. The quotient is therefore parametrised by the $A_m$-Weyl chamber. 
  
For simplicity, we assume the $m$ singular values are nonzero and $\mathbb{F} = \mathbb{R}$, and for cleaner narrative, we focus on the real case.

Since $K$ acts by isometry, we may compute the orbit mean curvature conveniently at \[M = L :=  \begin{bmatrix}
	\diag(\sigma_1, \cdots, \sigma_m) \\
	\mathbf{0}
\end{bmatrix} \] with and $\sigma_1>\cdots>\sigma_m > 0$. The vertical space at $M$, $\ver_M\mathbb{R}^{nm} = T_L K\cdot L$ is spanned by $\rho^L_*(A,0) = AL$ for $A\in \mathfrak{so}(n)$ and $\rho_*(0, A) = LA$ for $A\in \mathfrak{so}(m)$. 	We choose the basis \[\{(A_{ij}, \pm A_{ij})|i<j\leq m\}\cup\{(A_{ij}, 0)|i\leq m < j\}\] for the $K$-invariant field on $K/T$, $T$ being the maximal torus. Under $\rho^L_*$, \begin{equation}
(A_{ij}, A_{ij})\mapsto (\sigma_i+\sigma_j)A_{ij}, \text{ and } (A_{ij}, -A_{ij})\mapsto -(\sigma_i - \sigma_j)S_{ij}
\end{equation} for $i<j\leq m$, and 
\begin{equation}
(A_{ij}, 0) \mapsto -\sigma_i E_{ji}
\end{equation} for $i\leq m < j$. The horizontal space therefore consists of those of the form \[\begin{bmatrix}
D \\
\mathbf{0}
\end{bmatrix}, \] where $D>0$ is diagonal.

The above proves:
\begin{proposition}
	The mapping $\phi:\Mnmo \to \Mnmo/K$ is a Riemannian submersion. The quotient is parametrised by a single chart consisting of the $A_m$-Weyl chamber equipped with the Euclidean metric. 
\end{proposition}

Let $\lambda_1 = \sigma_1^2 > \cdots > \lambda_m > \lambda_{m+1} = 0 = \cdots = \lambda_n$ be eigenvalues of $WW^*$. 

\begin{proposition}
	The orbit volume function with respect to the quotient $\phi$ is 
	\[vol(K\cdot\Sigma) \propto \prod_{i\leq m; \; i < j}(\lambda_i - \lambda_j)^{2\beta}\] with $\beta = 1, 2, 3$ for real, complex, and quaternion respectively.
\end{proposition}

\begin{proof}

	For $\mathbb{F} = \mathbb{R}$ as shown above, the metric pulled back from the submanifold orbit $K\cdot\Sigma \subset (\Mnm, Fr)$ is diagonal with entries $(\sigma_i + \sigma_j)^2$ and $(\sigma_i - \sigma_j)^2$ for $i<j\leq m$, and $\sigma_j^2$ for $j>m$. The volume function is \begin{equation}
		vol(K\cdot\Sigma) = \det\phi^*g \propto \prod_{i\leq m; \; i < j}(\lambda_i - \lambda_j)^2.
	\end{equation}

	For $\mathbb{F} = \mathbb{C}$ or $\mathbb{H}$, we consider a basis for the imaginary parts
	\[\{(S_{ij}, \pm S_{ij})\}\] and observe that under $\rho^L_*$, \begin{equation}
		(S_{ij}, S_{ij}) \mapsto (S_{ij}L + L S_{ij}) = (\sigma_i + \sigma_j)S_{ij}
	\end{equation} and 
	\begin{equation}
		(S_{ij}, -S_{ij}) \mapsto (S_{ij}L - L S_{ij}) = (\sigma_i + \sigma_j)A_{ij},
	\end{equation} contributing to another factor of $(\lambda_i^2 - \lambda_j^2)$.
\end{proof}

\begin{corollary}
	The image process of $W$ through $\phi$ is \[d\sigma_i = d\omega_i +\beta\sum_{i\leq m; \; i < j}\frac{\sigma_i}{\lambda_i-\lambda_j}dt,\] with $\beta =1, 2, 4$ for real, complex, and quaternion respectively, where $\omega_i$'s are independent standard BM.
\end{corollary}

\begin{proof}
	The find the gradient, we compute partial derivatives  \begin{equation}
		\frac{\partial}{\partial\sigma_i} \log\prod_{i \neq j}(\lambda_i - \lambda_j)^{2\beta} = \sum_{j \neq i}\frac{2\beta\sigma_i}{\lambda_i-\lambda_j}.
	\end{equation}
\end{proof}

\begin{corollary}
	The eigenvalue processes of $W_tW_t^*$ satisfy \[ d\lambda_i = 2\sqrt{\lambda_i}d\omega_i + (\sum_{j\neq i}2\beta\frac{\lambda_i + \lambda_j}{\lambda_i - \lambda_j} + n)dt,\] recovering the formula derived in \cite{bru1989diffusions}.
\end{corollary}

\begin{proof}
	By the previous corollary,
\begin{align}
d\lambda_i & = d\sigma_i^2 = 2\sigma_id\sigma_i + d\sigma_id\sigma_i \\
& = 2\sigma_id\omega_i + \sum_{j\neq i}\frac{2\lambda_i}{\lambda_i - \lambda_j}dt  + 1dt \\
& = 2\sqrt{\lambda_i}d\omega_i + \sum_{j\neq i}(\frac{\lambda_i + \lambda_j}{\lambda_i - \lambda_j} + 1) + 1dt \\
& = 2\sqrt{\lambda_i}d\beta_i + (\sum_{j\neq i}\frac{\lambda_i + \lambda_j}{\lambda_i - \lambda_j} + n)dt.
\end{align} 
\end{proof}

\subsection{Dynkin's Brownian motion}

Dynkin's Brownian motion is defined under a (left-)invariant metric on $G = GL(n, \mathbb{F})$, \begin{equation}
	\langle A,B\rangle_M = \rtr M^{-1}A(M^{-1}B)^*,
\end{equation} for $M\in G$. Hence the corresponding left-invariant BM is given by \begin{equation}
	dM = M\strat W, 
\end{equation} or $M = e^W$, where $W$ is an $n\times n$-dimensional standard Brownian motion. We want to find out the eigenvalue for the process of positive definite matrix $X = MM^*$.
 
Same as the Wishart process, we can consider the singular spectrum of $M$. The Lie group $K = U(n)\times U(n)$ acts on $G$ by $(U, V): M \mapsto UMV^*$. The actions are still isometric: \begin{align}
	\langle UMAV^*, UMBV^*\rangle_{UMV^*} & = \rtr (VM^{-1}U^*UMAV^*)(VM^{-1}U^*UMBV^*)^* \\
	& = \tr AB ^* = \langle MA, MB\rangle_M.
\end{align} We consider the mapping $\phi: \Mn \to \Mn/K$, where the quotient is parametrised by the $A_n$-Weyl chamber where the singular values are distinct. 

As before, we may assume $M = L:= \diag(\sigma_1, \cdots, \sigma_n) > 0$ is diagonal with distinct entries and write $\rho^L: K \to K\cdot L$, and we deal with the case $\mathbb{F} = \mathbb{R}$ for now, remarking the other two cases in the end.

Consider the basis with image under $\rho_*$ \begin{equation}
\mathcal{B} =	\{(A_{ij}, -\frac{\sigma_i}{\sigma_j}A_{ij})\mapsto \Sigma(\frac{\sigma_j}{\sigma_i} - \frac{\sigma_i}{\sigma_j})E_{ij})\}\cup\{(A_{ij}, -\frac{\sigma_j}{\sigma_i}A_{ij})\mapsto \Sigma(\frac{\sigma_j}{\sigma_i} - \frac{\sigma_i}{\sigma_j})E_{ji})\}.
\end{equation}

 This shows the horizontal space is spanned by diagonal matrices and has dimension $n$. The quotient space is parametrised by the $A_n$-Weyl chamber, whose inner product at a point $\Sigma$ is \begin{equation}
\langle \dot{\Sigma_1}, \dot{\Sigma_2}\rangle_\Sigma = \tr\Sigma^{-1}\dot{\Sigma_1}\dot{\Sigma_2}\Sigma^{-1} 
\end{equation}  for two tangent vectors $\Sigma_1$ and $\Sigma_2$.

We conclude here for the quotient geometry:
\begin{proposition}\label{Dynkin geo}
	The mapping $\phi: GL(n) \to GL(n)/(U(n)\times U(n))$ is a Riemannian submersion, where $GL(n)$ is equipped with the (usual) left invariant metric. 
	
	The quotient space is parametrised by the $A_n$-Weyl chamber $\mathcal{W}$ with a single chart, on which the metric at a point  $\Sigma\in\mathcal{W}$ is \begin{equation*}
		\langle \dot{\Sigma_1}, \dot{\Sigma_2}\rangle_\Sigma = \tr\Sigma^{-1}\dot{\Sigma_1}\dot{\Sigma_2}\Sigma^{-1}. 
	\end{equation*} 
\end{proposition}

\begin{proposition}
	The orbit volume function with respect to the quotient $\phi$ satisfies 
	\[vol(K\cdot \Sigma) \propto \prod_{i \neq j}\frac{(\lambda_i - \lambda_j)^2}{\lambda_i\lambda_j}.\]
\end{proposition}

\begin{proof}
		For $\mathbb{F} = \mathbb{R}$, since the pull-back metric on $K/T$ is diagonal under the basis $\mathcal{B}$, \begin{equation}
		vol(K\cdot\Sigma) \propto \prod_{i \neq j} (\frac{\sigma_j}{\sigma_i} - \frac{\sigma_i}{\sigma_j})^2 = \prod_{i \neq j}\frac{(\lambda_i - \lambda_j)^2}{\lambda_i\lambda_j}.
	\end{equation}

	For $\mathbb{F} = \mathbb{C}$ or $\mathbb{H}$, similar to the Wishart case, we consider $S_{ij}$ in place of $A_{ij}$ for the imaginary parts.
\end{proof}

\begin{corollary}
Singular values of $M_t$ satisfy  \begin{equation*}
		d\sigma_i = \sigma_id\beta_i + \frac{1}{2}\sum_{j\neq i}((\frac{\lambda_i + \lambda_j}{\lambda_i - \lambda_j})^2\sigma_i - \frac{2(\lambda_i+\lambda_j)\lambda_j\sigma_i}{(\lambda_i-\lambda_j)^2})dt.
	\end{equation*}
\end{corollary}

\begin{proof}
	First we compute partial derivatives, \begin{align}
		\frac{\partial}{\partial\sigma_i}\log\prod_{i \neq j}\frac{(\lambda_i - \lambda_j)^2}{\lambda_i\lambda_j} & = \frac{\partial}{\partial\sigma_i}\sum_{1\leq i,j \leq n, i \neq j}(\log(\lambda_i-\lambda_j)^2 - \log \lambda_i\lambda_j) \\
		& = 2\sigma_i\sum_{j \neq i}(\frac{2}{\lambda_i-\lambda_j} - \frac{1}{\lambda_i}) \\
		& = \frac{4\sigma_i}{\lambda_i}\sum_{j \neq i}\frac{\lambda_i+\lambda_j}{\lambda_i-\lambda_j}.
	\end{align}
	
	Under the metric at $L$, the gradient is \begin{equation}
		\grad \log vol(K\cdot L) = \sum_i 4\sigma_i \sum_{j \neq i}\frac{\lambda_i+\lambda_j}{\lambda_i-\lambda_j}\partial_{\sigma_i}.
	\end{equation} 
	
\end{proof}

\begin{corollary}
	The eigenvalue processes of $M_tM_t^*$ satisfies \[d\lambda_i = 2\lambda_id\beta_i  + \sum_{j\neq i}\frac{(\lambda_i+\lambda_j)\lambda_i}{(\lambda_i-\lambda_j)^2}dt +\lambda_idt.\]
\end{corollary}

\begin{proof}

	\begin{align}
		d\lambda_i & = 2\sigma_id\sigma_i + d\sigma_id\sigma_i \\
				   & = 2\lambda_id\beta_i  + \sum_{j\neq i}((\frac{\lambda_i + \lambda_j}{\lambda_i - \lambda_j})^2\sigma_i - \frac{2(\lambda_i+\lambda_j)\lambda_j\lambda_i}{(\lambda_i-\lambda_j)^2})dt +\lambda_idt \\
				   & = 2\lambda_id\beta_i  + \sum_{j\neq i}\frac{(\lambda_i+\lambda_j)\lambda_i(\lambda_i+\lambda_j-2\lambda_j)}{(\lambda_i-\lambda_j)^2}dt +\lambda_idt \\
				   & = 2\lambda_id\beta_i  + \sum_{j\neq i}\frac{(\lambda_i+\lambda_j)\lambda_i}{(\lambda_i-\lambda_j)^2}dt +\lambda_idt.
	\end{align}
\end{proof}

\begin{remark}
	The form of the equation in  \cite{rogers_williams_2000} is given as
	\begin{equation}
		d\gamma_i = d\beta_i + \frac{1}{2}\sum_{j\neq i}\coth(\gamma_i-\gamma_j)dt =  d\beta_i + \frac{1}{2}\sum_{j\neq i}\frac{\lambda_i+\lambda_j}{\lambda_i-\lambda_j}dt,
	\end{equation} where $\gamma_i = \frac{1}{2}\log\lambda_i$ or $\lambda_i = e^{2\gamma_i}$. 
\end{remark}

\section{Eigenvector processes}\label{EV section}

In this section, we see how the eigenvectors evolve for the three processes described. We shall see that eigenvectors satisfy the SDE for BM on the fibre and are often not autonomous. Generic fibres are diffeomorphic to the group $K$ quotient by the maximal torus $T$, but the canonical diffeomorphisms do not ensure isometries. In other words, the evolution of eigenvectors involves eigenvalues and is not autonomous. 

Moreover, eigenvectors are not unique. The evolution depends on given initial solutions. With more geometric flavour, the initial solution is a ``lift'' to the group $K$.  


\subsection{Riemannian semidirect products}

First, let us reinvent some wheels with calculations on the local model of Riemannian submersions, which are locally trivial. 

\begin{definition}
	A product of two smooth manifolds $Y\times X$ is called a \textbf{Riemannian semidirect product}, denoted $Y \rtimes X$ if it is a Riemannian manifold such that the metric on the tangent space $T_{(y,x)}Y\oplus X$ is decomposed as $g_{x} + h$, where $g$ and $h$ are supported on $T_yY$ and $T_xX$ respectively, and $h$ is independent of $y\in Y$. 
\end{definition}

 Riemannian semidirect products are also a generalisation of  \textit{warped products}. What we have in mind is considering the $Y$ component as eigenvectors and $X$ as eigenvalues in the Riemannian submersions from the previous section. The name \textit{semidirect product} is used before the author finds a even more suitable one. It is not strictly related to the semidirect products of groups but only has similarity in ideas. 

\begin{proposition}\label{semidirect gamma}
	Given a Riemannian semidirect product as in the definition above. Using indices $a, b, c$ for coordinate functions of $Y$ and $i, j, k$ for those of $X$, we have Christoffel symbols as below.
	
	The Levi-Civita connection for vectors along $X$ has components
	\begin{equation}
		\Gamma_{ij}^k = \tilde{\Gamma}_{ij}^k,
	\end{equation}
	and \begin{equation}
		\Gamma_{ij}^a = 0,
	\end{equation} where $\tilde{\Gamma}$ is the Christoffel symbol for $(X, h)$.

	The cross-term derivatives are \begin{equation}
		\Gamma_{ia}^b = \frac{1}{2}g^{bc}g_{ac, i}
	\end{equation} and 
	\begin{equation}
		\Gamma_{ia}^j = \frac{1}{2}h^{jk}h_{ki,a}.
	\end{equation}

	Finally along $Y$, 
	\begin{equation}
		\Gamma_{ab}^c = \hat{\Gamma}_{ab}^c,
	\end{equation} where $\hat{\Gamma}$ denotes the Christoffel symbol for $(Y, g_x)$, and 
	\begin{equation}
		\Gamma_{ab}^i = -\frac{1}{2}h^{ij}g_{ab,i}.
	\end{equation}
\end{proposition}

\begin{proof}
	These come directly from the well-known equation
	\begin{equation}
		\Gamma_{ij}^k = \frac{1}{2}g^{kl}(g_{li,j} + g_{jl,i} - g_{ij,l}), 
	\end{equation} noting that the metric is of the form
	\[\begin{bmatrix}
		g & 0 \\
		0 & h
	\end{bmatrix},\] where there are no cross terms.
\end{proof}

\begin{corollary}\label{prX prY}
	Given the Riemannian semidirect product as above. Let $\pr_Y$ and $\pr_X$ be projections to each component respectively. Then $\pr_Y$ has totally geodesic fibres, and $\pr_X$ is a Riemannian submersion. 
\end{corollary}

\begin{proof}
	The first two equations from the proposition above show $\{y\}\times X$ has nonvanishing second fundamental form. Since $h$ does not depend on $y\in Y$, $\pr_X$ is a Riemannian submersion. 
\end{proof}

The following theorem indicates the last element we need is the metric tensor of each fibre.

\begin{proposition}\label{semidirect BM}
	The SDE for Brownian motion on a Riemannian semidirect product as notated above is a sum of three terms: (1) BM on $(Y, g)$ that is dependent on $x\in X$, (2) BM on $(X,h)$, and (3) a drift term proportional to the mean curvature of $Y\times \{x\} \subset Y\rtimes X$.
\end{proposition}

\begin{proof}
	With the formula for BM in local coordinates (Cf. \cite{hsu2002stochastic}), we see 
	\begin{equation}
		dy^c = (\sqrt{g})^c_a|_{(y,x)} dB^a - \frac{1}{2}g^{ab}\Gamma_{ab}^c dt,
 	\end{equation} which is BM on $(Y, g_x )$, and
 	\begin{equation}
 		dx^k = (\sqrt{h})^k_i(x) dB^i -\frac{1}{2}h^{ij}\Gamma_{ij}^k dt - \frac{1}{2}g^{ab}\Gamma_{ab}^k dt.
 	\end{equation} The first two terms describe BM on $(X, h)$ whereas the last term is a multiple of the mean curvature of $Y\times\{x\}$ at $(y,x)$.
\end{proof}

\begin{corollary}
	Given the Riemannian semidirect product as previously described and a BM on it denoted $M_t$. Then at time $t$, the ``vertical process'' $\pr_Y(M_t)$ satisfies the SDE for BM on the fibre $\pr_X^{-1}\circ\pr_X(M_t)$.
\end{corollary}

\begin{proof}
	By Cor. \ref{prX prY}, $\pr_Y$ has totally geodesic fibres. The statement follows combined with Prop. \ref{semidirect BM}.
\end{proof}

Taking Dyson BM for example, a fibre of the submersion $\phi: \hermo \to \mathcal{W}$ is $K=U(n)$ quotient by its maximal torus $T$. But what we want is a motion on $K$. This requires a lift $K\to K/T$ and eventually a lift $K\times \mathcal{U} \to K/T\times \mathcal{U}$ for an open set $\mathcal{U}\subset\mathcal{W}$ over which the submersion $\phi$ is trivial, which are Riemannian submersions. The following fact ensures the existence of valid lifts, also noting that we shall see the metrics on fibres are still $K$-invariant.  

\begin{proposition}\label{reductive lift}(C.f. \cite{gallier2020differential})
	
	Let the space of right cosets, $G/K$, be a reductive homogeneous space as described in the definition above. If $\mathfrak{m}$ has an $Ad(K)$-invariant inner product $\langle\cdot,\cdot\rangle_\mathfrak{m}$, then is quotient map $G \to G/K$ is a Riemannian submersion where $G$ is endowed with any left $G$-invariant metric $\langle\cdot, \cdot\rangle_\mathfrak{g}$ extending $\langle\cdot,\cdot\rangle_\mathfrak{m}$ such that $\mathfrak{k}$ is orthogonal to $\mathfrak{m}$, where $\mathfrak{m}$ maps isometrically to $T_K G/K$.

\end{proposition}

 Note that, taking Dyson BM for example, when $\mathbb{F} = \mathbb{R}$, $O(n) \to O(n)/T$ is a covering map since $T$ is zero-dimensional. For $\mathbb{C}$ and $\mathbb{H}$, the lifted metric for $U(n)\to U(n)/T$ is however not unique in the fibre direction. Potentially sub-Riemannian structure might be a more natural tool.


\subsection{Eigenvecors of Dyson Brownian motion}\label{Dyson EV}

We stick to the complex case for Dyson BM. The other two cases are similar. 

The only hands-on work required is describing the metric on each fibre by writing down an orthonormal basis of $\mathfrak{su}(n)$.

Define 
\begin{equation}
	\tilde{A}_{ij} = \frac{E_{ij} - E_{ij}}{\sqrt{2}(\lambda_j-\lambda_i)},
\end{equation} 
\begin{equation}
	\tilde{S}_{ij} = \frac{E_{ij} + E_{ji}}{\sqrt{2}(\lambda_j-\lambda_i)},
\end{equation} for $i\neq j$. Then we have that
\[\{[\tilde{A}_{ij},\Lambda]  = \tilde{S}_{ij}\}\cup\{[\ib\tilde{S}_{ij}, \Lambda] = \ib A_{ij}\}(\cup\{[\jb\tilde{S}_{ij}, \Lambda] = \jb A_{ij}\}\cup\{[\kb\tilde{S}_{ij}, \Lambda] = \kb A_{ij}\})\] for $i<j$ is an o.n.b. for $\ver_\Lambda\hermo$. 

We may endow the canonical metric for the maximal torus $T< K$ and the horizontal space to be the spanned of the basis above. This is a valid lift since the maximal torus commutes with diagonal matrices and conditions in \ref{reductive lift} are satisfied. BM on the Lie algebra under such metric is thus $A:= \tilde{A}_{ij}\strat W_1^{ij} + \ib\tilde{S}_{ij}\strat W_2^{ij} + \ib(D_i-D_j)\strat W_3^{ij}$, where $W$'s are independent standard BM. Brownian motion on $K$ under such metric satisfies the following SDE.
\begin{align}
	dQ & = Q \strat A \\
	& = QdA + \frac{1}{2}dQdA \\
	& = QdA + \frac{1}{2}QdAdA \\
	& = QdA + \frac{1}{2}Q\diag(-\sum_{j\neq 1}\frac{dt}{2(\lambda_1-\lambda_j)^2}, \cdots) \\
	& =  QdA - \frac{1}{4}Q\diag(\sum_{j\neq 1}\frac{dt}{(\lambda_1-\lambda_j)^2}, \cdots) \\
\end{align} for $\mathbb{F} = \mathbb{R}$. 

Let $q_i$ be the $i$-th column of $Q$, we have the evolution of a single eigenvector.

\begin{proposition} 
	An eigenvector of a real Dyson BM satisfies 
	\[dq_i = \frac{1}{\sqrt{2}}\sum_{j\neq i} \frac{1}{\lambda_i-\lambda_j}\cdot q_i dW^{ij} - \frac{1}{4}\sum_{j \neq i}\frac{1}{(\lambda_i- \lambda_j)^2}\cdot q_i dt.\]
\end{proposition}

The complex and quaternion cases are similar.

\subsection{Eigenvectors of Wishart process}\label{Wishart EV}

Observe that
 \begin{equation}
	\frac{1}{\sigma_i+\sigma_j}\rho^L_*(A_{ij}, A_{ij}^*) = \frac{A_{ij}L + LA_{ij}}{\sigma_i + \sigma_j} = A_{ij}
\end{equation} and
\begin{equation}
	\frac{1}{\sigma_i-\sigma_j}\rho^L_*(A_{ij}, A_{ij}) = \frac{A_{ij}L - LA_{ij}}{\sigma_i - \sigma_j} = S_{ij}
\end{equation} for $i <j\leq m$, and 
\begin{equation}
	\rho^L(\frac{A_{ij}}{\sigma_i}, 0) = -E_{ji}
\end{equation} for $i\leq m < j$. They are orthogonal under the Frobenius norm and span the vertical space. (Since $\rho^L(A_{ij}, 0) = 0$.) So we have an o.n.b.  $\{\frac{1}{\sqrt{2}(\sigma_i+\sigma_j)}(A_{ij}, A_{ij}^*)\}\cup\{\frac{1}{\sqrt{2}(\sigma_i-\sigma_j)}(A_{ij}, A_{ij})\}$.

We are only interested in the eigenvectors of $WW^*$, which comes from the corresponding Lie algebra BM on $\mathfrak{su}(n)$ \begin{align}
	A & = W_1^{ij}\frac{A_{ij}}{\sqrt{2}(\sigma_i-\sigma_j)} + W_2^{ij}\frac{A_{ij}}{\sqrt{2}(\sigma_i-\sigma_j)} \\
	& = \frac{(\sigma_j - \sigma_i)W_1^{ij} + (\sigma_i + \sigma_j)W_2^{ij}}{\sqrt{2}(\lambda_j - \lambda_i)}A_{ij} \\
	& =: \frac{\sqrt{(\sigma_j - \sigma_i)^2 + (\sigma_i + \sigma_j)^2} }{\sqrt{2}(\lambda_j - \lambda_i)} Z^{ij}A_{ij} \\
	& = \frac{\sqrt{\lambda_i + \lambda_j}}{(\lambda_j - \lambda_i)}Z^{ij}A_{ij},
\end{align} where $Z$'s are independent standard BM after a linear transform of $W$'s.

The eigenvector process from $WW^*$ is thus \begin{align}
	dQ & = Q\strat A \\
	& = QdA + \frac{1}{2}QdAdA\\
	& = QdA - \frac{1}{2}Q\diag(\sum_{j \neq 1}\frac{\lambda_1 + \lambda_j}{(\lambda_1 - \lambda_j)^2}, \cdots)dt.
\end{align} It recovers the result from \cite{bru1989diffusions} as follows.

\begin{proposition}
	The $i$-th eigenvector , $q_i$, of Wishart's process satisfies 
	\[dq_i = \sum_{j \neq i}\frac{\sqrt{\lambda_i + \lambda_j}}{(\lambda_j - \lambda_i)}\cdot q_jdZ^{ij} - \frac{1}{2}\sum_{j\neq i}\frac{(\lambda_i + \lambda_j)}{\lambda_i - \lambda_j}\cdot q_idt\]
\end{proposition}

\subsection{Eigenvecors of Dynkin's Brownian motion}\label{Dynkin EV}
We adjust the notation by denoting $A_{ij} = E_{ij} - E_{ji}$ and see that  
\begin{equation}
	\rho^L_*: \tilde{A}_{ij} := (\frac{\sigma_i\sigma_j}{\lambda_j-\lambda_i}A_{ij}, - \frac{\lambda_i}{\lambda_j-\lambda_i}A_{ij}^* ) \mapsto \frac{\sigma_i\sigma_j}{\lambda_j-\lambda_i}A_{ij}L - \frac{\lambda_i}{\lambda_j-\lambda_i}LA_{ij} = LE_{ij}
\end{equation} and \begin{equation}
	\rho^L_*:  \hat{A}_{ij} :=  (\frac{\sigma_i\sigma_j}{\lambda_j-\lambda_i}A_{ij}, - \frac{\lambda_j}{\lambda_j-\lambda_i}A_{ij}^* ) \mapsto \frac{\sigma_i\sigma_j}{\lambda_j-\lambda_i}A_{ij}L - \frac{\lambda_j}{\lambda_j-\lambda_i}LA_{ij} = LE_{ji}.
\end{equation} Therefore $\{\tilde{A}_{ij}\}\cup\{\hat{A}_{ij}\}$ is an orthonormal basis under the geometry described in \ref{Dynkin geo}.

BM on $\mathfrak{k} = Lie(O(n)\times O(n))$ is therefore \begin{align}
	\tilde{A}_{ij}W_1^{ij} + \hat{A}_{ij}W_2^{ij} & = (\frac{\sigma_i\sigma_j}{\lambda_i - \lambda_j}(W_1^{ij} + W_2^{ij})A_{ij}, \frac{\lambda_iW_1^{ij + \lambda_jW_2^{ij}}}{\lambda_j - \lambda_i}A_{ij}) \\
	& =: (\frac{\sqrt{2}\sigma_i\sigma_j}{\lambda_j - \lambda_i}Z_1^{ij}A_{ij}, \frac{\sqrt{\lambda_i^2 + \lambda_j^2}}{\lambda_j - \lambda_i}Z_2^{ij}A_{ij})\\
	& =:(A,B) ,
\end{align} where $W$'s are independent standard BM, and $Z$'s are independent standard BM from a linear transformation of $W$'s as in the equations. 

The SDE for the vertical process $(U,V)$ goes as 
\begin{align}
	dU & = U\strat A = U dA + \frac{1}{2}UdAdA \\
	 & = U dA - U\diag(\cdots, \sum_{j \neq i}\frac{\lambda_i \lambda_j}{(\lambda_i - \lambda_j)^2}, \cdots)dt, 
\end{align} which corresponds to the process $N$ in (9.2) of \cite{article}, and 

\begin{align}
	dV & = V\strat B = V dB + \frac{1}{2}VdBdB \\
		& = V dB - \frac{1}{2}V\diag(\cdots, \sum_{j \neq i}\frac{\lambda_i^2 + \lambda_j^2}{(\lambda_i - \lambda_j)^2}, \cdots) dt, 
\end{align} which is the process $M$ in (6.13) or (8.1) in \cite{article}, noting that (6.11) and (6.12) in the paper ensure $\tilde{S}$ and $\tilde{A}$ are ''standard'' symmetric and skew-symmetric BM. Moreover, since \begin{equation}
\sqrt{ 1 + (\frac{\lambda_i + \lambda_j}{\lambda_j - \lambda_i})^2} = \frac{\sqrt{\lambda_i^2 + \lambda_j^2}}{\lambda_j - \lambda_i}, 
\end{equation} it translates to the form we derived.

\section{Constructions for general $\beta$-processes}\label{general beta section}

In our previous paper we constructed the matrix model whose eigenvalue process is the $\beta$-Dyson BM for $\beta>0$. Here we see how it works in general.

\begin{lemma}
	Let $\phi: M \to N$ be a submersion of smooth manifolds. If $X$ is a process on $M$ solves an SDE \[dM = V_i\strat Z^i\] for some real semimartingales $Z_i$ and vector fields $V_i$, then the image process $\phi(X)$ satisfies \[d\phi(X) = \phi_*(V_i)\strat Z^i.\]
\end{lemma}

\begin{proof}
	Let $f$ be an arbitrary smooth function on $N$. Since $\phi$ is surjective, $f\circ\phi$ is a smooth function on $M$ constant along each fibre. Thus \begin{equation}
		df\circ\phi(X) = V_i f\circ\phi(X)\strat Z^i = \phi_*(V_i)f(X) \strat Z^i.
	\end{equation}
\end{proof}

With respect to a Riemannian submersion $\phi: M\to N$, we introduce the operator $\prv$ on $TM$ to be the orthogonal projection to the vertical distribution and $\prh$ the orthogonal projection to the horizontal distribution. To ``project'' a stochastic differential $dX$ on $M$, Nash embedding allows us to assume $M$ is a submanifold of Euclidean space, and the projections are realised as projection matrices operating on the vector $dX$. Equivalently 

\begin{proposition} 
	Let $\phi: M^{n+l} \twoheadrightarrow N^n$ be a submersion of smooth manifolds that is a local trivial fibration with relative dimension $l$. Suppose furthermore $\phi$ is a Riemannian submersion. Let $W_M$ be a Brownian motion on $M$, and let $\prv$ and $\prh$ be orthogonal projection to the vertical and the horizontal directions respectively. Define \[dX = \alpha\prv dW_M + \beta\prh dW_M\] for $\alpha, \beta \in \mathbb{R}$. Then the image process satisfies \[d\phi(X) = \alpha \frac{l}{2}\phi_*(\vec{H}) + \beta dW_N, \] where $\vec{H}$ is the fibre mean curvature and $W_N$ is a BM on $N$.
\end{proposition}

\begin{proof}
	We may assume $M$ is isometrically embedded in Euclidean space $\mathbb{R}^s$ with standard basis $\{e_i\}_{i = 1, \cdots, s}$ and the origin $o\in M$ such that  $e_1, \cdots, e_{n+l}$ span $ToM$ and $e_1, \cdots, e_l$ spans $\ver_oM$. 
	
	Let $W$ be a Euclidean BM on $\mathbb{R}^s$, $\xi_i$ be the orthogonal projection of $e_i$ to $TM$, and $\zeta_i = \prv e_i$. Note that $\prv \xi_i = \zeta_i$ since $\ver M\subset TM$. Then the solution of $dY = \xi_i\strat W^i$ is a BM on $M$. 
	
	Let \begin{equation}
		dX = \zeta_i d Y^i = \zeta_i \strat W^i - \frac{1}{2}\nabla_{\zeta_i}\zeta_jj dY^idY^j.
	\end{equation} So \begin{equation}
	d\phi(X) = -\frac{1}{2}\sum_i\phi_*(\nabla_{\zeta_i}\zeta_j) dt = -\frac{l}{2}\phi_*(H)
	\end{equation}

	This further shows that if $d\tilde{X} = \prh dY$, then $\phi(\tilde{X})$ is a BM on $N$. The proposition follows by the additivity of integrals. 
\end{proof} 

\begin{remark}
	The proposition extends Lemma 3.1 in \cite{https://doi.org/10.48550/arxiv.2205.06737} and corrects its proof.
\end{remark}

\begin{remark}(Construction via metric change)
	Decompose the metric $g$ into vertical and horizontal parts $g = g_h + g_v$ and define a new metric $\tilde{g} = r^2 g_h + g_v$. Let $\tilde{W}$ be the new BM under $\tilde{g}$. Then \begin{equation}
		d\phi(\tilde{W}) = rdW_N - \frac{l}{2r}\phi_*(H)
	\end{equation} since fibre volume function stays the same. 
	
\end{remark}


\section{Concluding remarks}
We compare the  three eigenvalue processes in the table below.
\vspace{0.5cm}

{\centering
	\begin{tabular}{|c||c|c|c|}
		\hline 
		& Dyson & Wishart & Dynkin \\
		\hline \hline
		ensemble & Gaussian & Laguerre & $e^W e^{W^*}$ \\
		\hline
		domain & hermitian & nonnegative definite & positive definite \\
		\hline
		geometry & Euclidean & Bures-Wasserstein & Cartan-Hadamard \\
		\hline
		invariant group &  $U(n)$ & $U(n)\times U(m)$ & $U(n)\times U(n)$ \\
		\hline
		orbit volume &　$\prod_{i \neq j}(\lambda_i-\lambda_j)^\beta$ & $\prod_{i \neq j} (\lambda_i - \lambda_j)^{2\beta}$ & $\prod_{i \neq j} (\frac{\lambda_i-\lambda_j}{\sqrt{\lambda_i\lambda_j}})^{2\beta}$ \\
		\hline
		
	\end{tabular}
} 
\vspace{0.5cm} 	

\begin{itemize}
	\item The first row shows the ensembles when we take the time at $t=1$ of the dynamical versions. The ensemble from Dynkin's BM is by taking $W$ to be an $n\times n$ $\mathbb{F}$-standard normal. However, I have not found further study of the ensemble in literature, therefore the lack of a name.
	
	
	\item The \textit{geometry} row shows the related Riemannian structures under which the matrix processes are BM. 
	
	We have seen clearly Dyson's hermitian process is the BM on the linear subspace of hermitian matrices. 
	
	The Wishart process is not exactly the BM under the Bures-Wasserstein geometry, which is defined as the quotient geometry from the Frobenius geometry by the unitary group, whose orbits have nonzero mean curvatures, and hence by \ref{BM projection formula}, the mean curvature term shows up in the image process. See \cite{https://doi.org/10.48550/arxiv.2205.06737} for more details. 
	
	The Cartan-Hadamard geometry is defined as the quotient of the usual left invariant metric of $GL(n)$ by the unitary group, which is also mentioned in \cite{https://doi.org/10.48550/arxiv.2205.06737}.
	
	These three are the most common and natural Riemannian structure on positive definite matrices, with flat, nonnegative \cite{massart2019curvature}, and negative \cite{lang2012fundamentals} curvature respectively. 
\end{itemize}

\subsection{Future directions}

\subsubsection{Towards universality?}
Geometry is but an illusion built on the basis of analysis and algebra. However the intuitions coming from it often help greatly on our cognitive process, as what has happened so far in this paper. With optimism, we could hope this to be only the beginning of what geometry can help understanding RMT. 

We are of course curious what geometric meanings other objects and concepts in RMT have, eventually aiming to shed lights on central questions such as universality conjecture. For instance, one may solve the Fokker-Planck equation for the Dyson BM \ref{Dyson BM formula} by factoring out the Vandermonde $\prod_{i \neq j}(\lambda_i-\lambda_j)$ to reduce to the heat equation. Now we know the Vandermonde is proportional to the orbit volume \ref{Dyson orbit volume}, do we have more geometric interpretations? 


\subsubsection{What happens at singularity?}

Eigenvalues collide when the orbit is degenerate. Although this is but an interpretation of eigenspace and not completely new to us, the picture is inviting for geometers to think of resolution of singularities. Roughly, for submersions of the form $\phi: \mathcal{M} \to \mathcal{W}$, mapping from matrices to spectra, we consider (at least locally) surjective mappings: 

\begin{equation}
	K/T\times \mathcal{U}  \xrightarrow{\xi} \mathcal{R} \xrightarrow{\phi} \phi(\mathcal{R}), 
\end{equation} for open sets $\mathcal{U} \subset \mathcal{W}$ and $\mathcal{R} \mathcal{M}$, noting that $K/T$ is a generic fibre ($T$ being the maximal torus). 

There are at least a couple reasons this is of interests. We are curious to see what happens at singularities. \cite{allez2013diffusive} For designing algorithms, it might be desirable to have some ways to pass through singularities.

\subsubsection{Mean curvature control system}
In \cite{https://doi.org/10.48550/arxiv.2205.06737}, we use the metric as the control on the mean curvature flow resulting from projecting BM through a submersion. More specifically, we consider Riemannian submersions $\phi_c: (M,g_c) \to (N,h_c)$ with the control $c$ from some set $\mathcal{C}$, resulting in orbit mean vector field $\vec{H}_c$, and hence the image of BM $W_t \in (M,g_{c(t)})$ drifts with (up to a scale) $\phi_*(\vec{H}_{c(t)}(W_t))$. Thus we can drive the eigenvalues by varying the metric of the total space consisting of matrices. 

\appendix

\section{Derivations of mean curvature according to definition}

We include the computations for orbit mean curvature directly from definition. 

\subsection{Dyson}

\begin{lemma}
	The orbit mean curvature for Dyson Brownian motion is \[\vec{H}(U\Lambda U^*) = U\diag(\sum_{j \neq 1}\frac{\beta}{\lambda_1-\lambda_j}, \cdots)U^*.\]
\end{lemma}

\begin{proof}
	Let $U_t = \exp(tA)$ for $A \in \algsu$ and $S_t:= U_t\Lambda U_t^*$. Then $S_t$ is a trajectory of the vector field $\phi_*(UA) = [A, S_t]$. Therefore \begin{equation}
		\nabla_{\phi_*(U_tA)}\phi_*(U_tA)|_\Lambda = \frac{d^2S_t}{dt^2}|_0 = [A, [A,\Lambda]].
	\end{equation} 	

	Using the o.n.b. in \ref{Dyson EV},  
	\begin{equation}
		[\tilde{A}_{ab},[\tilde{A}_{ab}, \Lambda]] = \frac{D_a - D_b}{\lambda_b-\lambda_a} \in \ver_\Lambda \hermo
	\end{equation} and 
	\begin{equation}
		[\ib\tilde{S}_{ab}, [\ib\tilde{S}_{ab}, \Lambda]] = -\frac{D_b-D_a}{\lambda_b-\lambda_a} \in \ver_\Lambda \hermo
	\end{equation}  Let $c = \dim_{\mathbb{R}} U(n) - n$ be the real dimension of each orbit. The mean curvature is hence
	\begin{align}
		H = & \frac{1}{c}\sum_{a<b}[i\tilde{S}_{ab}, [i\tilde{S}_{ab}, \Lambda]] + \frac{1}{c}\sum_{a<b}[\tilde{A}_{ab}, [\tilde{A}_{ab}, \Lambda]] \\ 
		= & \frac{2}{c}\sum_{a<b} \frac{1}{\lambda_a-\lambda_b}(D_b-D_a) \\
		=  & -\frac{2}{c}\diag(\cdots, \sum_{b\neq a}\frac{1}{\lambda_a - \lambda_b}, \cdots)
	\end{align}
\end{proof}

\subsection{Wishart}

\begin{lemma}
	The orbit mean curvature for the Wishart process is \[\vec{H}(ULV^*) = U\diag(\sum_{j \neq 1}\frac{\beta}{\lambda_1-\lambda_j}, \cdots)V^*.\]
\end{lemma}

\begin{proof}
	
	Let $c = \dim U(n) + \dim U(m) - n$ and use the o.n.b. in \ref{Wishart EV}. Similarly to the previous section, the mean curvature
	\begin{equation}
		c\vec{H} = \sum_{i<j}\frac{A_{ij}^2L + 2 A_{ij}LA_{ij} + LA_{ij}^2}{(\sigma_i+\sigma_j)^2} 
		+ \sum_{i<j}\frac{A_{ij}^2L - 2 A_{ij}LA_{ij} + LA_{ij}^2}{(\sigma_i-\sigma_j)^2}.
	\end{equation}
	
	Note $A_{ij}^2L = LA_{ij}^2 = -\frac{1}{2}(\sigma_i D_i + \sigma_j D_j)$ and $A_{ij}LA_{ij} = -\frac{1}{2}(\sigma_j D_i + \sigma_i D_j)$. Hence the summand equals \begin{equation}
		-(\frac{1}{(\sigma_i+\sigma_j)^2} + \frac{1}{(\sigma_i-\sigma_j)^2})(\sigma_i D_i + \sigma_j D_j)
		-(\frac{1}{(\sigma_i+\sigma_j)^2} - \frac{1}{(\sigma_i-\sigma_j)^2})(\sigma_j D_i + \sigma_i D_j)
	\end{equation}
	\begin{align}
		= & -\frac{1}{(\lambda_i-\lambda_j)^2}(2(\lambda_i + \lambda_j)(\sigma_i D_i + \sigma_j D_j) - 4\sigma_i\sigma_j (\sigma_j D_i + \sigma_i D_j)) \\
		= & -\frac{2}{(\lambda_i-\lambda_j)^2}((2\lambda_j\sigma_i-\lambda_j\sigma_i-\lambda_i\sigma_j)D_i + (2\lambda_i\sigma_j - \lambda_j\sigma_i - \lambda_i\sigma_j)D_j) \\
		= & -\frac{2}{\lambda_i-\lambda_j}(\sigma_iD_i - \sigma_jD_j). 
	\end{align}
	
\end{proof}

\subsection{Dynkin}

\begin{lemma}
	The orbit mean curvature for Dynkin's BM is \[\vec{H}(ULV^*) = U\diag(\sum_{j\neq 1}((\frac{\lambda_1 + \lambda_j}{\lambda_1 - \lambda_j})^2\sigma_1 - \frac{2(\lambda_1 +\lambda_j)\lambda_j\sigma_1}{(\lambda_1 - \lambda_j)^2}), \cdots)V^*.\]
\end{lemma}

\begin{proof}

	Let $c = \dim O(n) - n$ and consider the o.n.b. in \ref{Dynkin EV}. Similarly, 
	\begin{align}
		c\vec{H} & = & \sum_{i<j}\frac{\lambda_i\lambda_j}{(\lambda_i-\lambda_j)^2}A_{ij}^2  - 2\frac{\sigma_i\sigma_j\lambda_i}{(\lambda_i-\lambda_j)^2}A_{ij}LA_{ij} + (\frac{\lambda_i}{\lambda_i-\lambda_j})^2LA_{ij}^2 \\
		& & + \sum_{i<j} \frac{\lambda_i\lambda_j}{(\lambda_i-\lambda_j)^2}A_{ij}^2  - 2\frac{\sigma_i\sigma_j\lambda_j}{(\lambda_i-\lambda_j)^2}A_{ij}LA_{ij} + (\frac{\lambda_j}{\lambda_i-\lambda_j})^2LA_{ij}^2 \\
		& = & -(\frac{\lambda_i + \lambda_j}{\lambda_i - \lambda_j})^2(\sigma_iD_i + \sigma_jD_j) - (\frac{2\sigma_i\sigma_j(\lambda_i+\lambda_j)}{\lambda_i - \lambda_j})^2(\sigma_jD_i + \sigma_iD_j) \\
		& = & -(\frac{\lambda_i + \lambda_j}{\lambda_i - \lambda_j})^2(\sigma_iD_i + \sigma_jD_j) - (\frac{2\sigma_i\sigma_j(\lambda_i+\lambda_j)}{\lambda_i - \lambda_j})^2(\lambda_j\sigma_iD_i + \lambda_i\sigma_jD_j).  
	\end{align}

\end{proof}

\bibliographystyle{alpha}
\nocite{*}
\bibliography{DysonRef}

\end{document}